\documentclass{conm-p-l}
\usepackage{amsmath}

\newtheorem{theorem}{Theorem}[section]
\newtheorem{proposition}[theorem]{Proposition}
\newtheorem{corollary}[theorem]{Corollary}
\newtheorem{cor}[theorem]{Corollary}

\newtheorem{notation}[theorem]{Notation}

\newtheorem{remark}[theorem]{Remark}
\newtheorem{lemma}[theorem]{Lemma}
\newtheorem{definition}[theorem]{Definition}

\newenvironment{rmk}{\begin{remark}\rm}{\end{remark}}
\newenvironment{defn}{\begin{definition}\rm}{\end{definition}}
\newenvironment{notn}{\begin{notation}\rm}{\end{notation}}

\newcommand{\R}{\mathbb{R}}
\newcommand{\Z}{\mathbb{Z}}
\newcommand{\N}{\mathbb{N}}
\newcommand{\Q}{\mathbb{Q}}

\newcommand{\e}{\mathbf{e}}

\newcommand{\lam}{\ensuremath{{\lambda}}}
\newcommand{\NP}{\mathrm{NP}}

\newcommand{\aff}{\mathrm{aff}}
\newcommand{\conv}{\mathrm{conv}}
\newcommand{\pos}{\mathrm{pos}}

\newcommand{\init}{\ensuremath{\mathrm{in}}}

\newcommand{\xvec}[1]{\ensuremath{x_{1}, \ldots, x_{#1}}}

\newcommand{\lcm}{\ensuremath{\mathrm{lcm}}}
\newcommand{\grp}{\ensuremath{\mathrm{grp}}}

\newcommand{\blam}{\boldsymbol{\lambda}}
\newcommand{\blambda}{\boldsymbol{\lambda}}
\newcommand{\bomega}{\boldsymbol{\omega}}
\newcommand{\bbeta}{\boldsymbol{\beta}}
\newcommand{\balpha}{\boldsymbol{\alpha}}
\newcommand{\bgamma}{\boldsymbol{\gamma}}

\begin{document}
 
\title{Some Normal Monomial Ideals}
\author{Marie A. Vitulli}
\address{Department of Mathematics\\1222 University of Oregon, Eugene, OR 
97403-1222}
\email{vitulli@math.uoregon.edu}
\subjclass{Primary 13C13; Secondary 13A20, 13F20}
\thanks{This article is dedicated to our friend and colleague, 
Ruth I. Michler, whose energy and presence enriched our community.}

\begin{abstract}
    In this paper we investigate the question of normality for 
    special monomial ideals in a polynomial ring over a field.  We first 
    include some expository sections that give the basics on the 
    integral closure of a ideal, the Rees algebra on an ideal, and 
    some fundamental results on the integral closure of a monomial 
    ideal.
    
\end{abstract}
\maketitle

\section{Introduction}\label{sec:intro}

There has always been considerable interest in monomial ideals in a polynomial
ring in $n$ indeterminates.  Many properties of the ideal $I$ can be translated 
into properties about the set $\Gamma(I)$  of exponents of monomials in 
the ideal and its convex hull $\NP(I)$ (the \emph{Newton polyhedron 
of $I$})
and these sets can largely be understood by the convex geometry 
of $\R^{n}$.  For example, it is well known that the monomial ideal $I$ is
integrally closed if and  only if $\Gamma(I) = \NP(I) \cap \N^{n}$
(see Theorem \ref{normalization}).   


Another algebraic notion that has significant geometric consequences 
is that of a normal ideal.  An ideal $I$ of a ring $R$ is said to be 
normal provided that all its positive powers are integrally closed. 
If $R$ is a normal integral domain then the ideal $I$ in $R$ is 
normal if and only if the Rees algebra $R[It]$ is normal.  The Rees
algebra is the algebraic counterpart to blowing up a scheme along a
closed subscheme so various geometric properties of the blow-up can be
interpreted as algebraic properties of the Rees algebra.   If
$I$ is a  monomial ideal and $R$ is a finite-dimensional polynomial ring
over a  field we can again translate the question of normality of $I$
into  questions about the exponent set $\Gamma(I)$ and the Newton 
polyhedron $\NP(I)$ of $I$.  However, no concise answer to the question
of  when $I$ is normal is known, nor does it appear likely that one can 
ever be given.  The question of normality of monomial ideals is the 
main focus of this paper.  We spend most of our time studying a special
class of monomial ideals that are primary to the ideal generated by the
indeterminates.

We now give a brief account of the organization of this manuscript. 
Our aim is to make this paper largely self-contained, thus we include 
proofs of many results that are well-known as ``folklore'' but whose 
proofs are elusive.   In Section
2.1 we present several preliminary results about the integral closure of graded
rings. We give results for both the relative and absolute notions of integral
closure.   In Section 2.2 we present results on the integral closure of a
homogeneous ideal in a graded ring. In Section 3 we specialize to the case  of
monomial ideals in a polynomial ring in $n$ indeterminates over a field $K$. 
By introducing an $\N^n$-grading on $K[\xvec{n}]$ we can identify monomial
ideals and homogeneous ideals of $K[\xvec{n}]$.  
We then present some well-known results about monomial ideals. We start 
with a complete account of the fact that the integral closure of  a
monomial ideal is again a monomial ideal (see Corollary
\ref{cor:monideal}) whose exponent set is the  set of integral points in
the Newton polyhedron of the ideal (see Theorem \ref{normalization}).  In
Section 3.1 we   turn our attention to the question of when a monomial
ideal is normal.   In section 3.2 we are interested in a special class of
monomial ideals that arise  from vectors $\blam = (\lam_{1}, \ldots ,
\lam_{n})$ of positive  integers.  First define $J(\blam) =
(x_{1}^{\lam_{1}}, \ldots,  x_{n}^{\lam_{n}})$ and then define
$I(\blam)$ to be the integral  closure of $J(\blam)$.
We introduce several recent results of Reid, Roberts, and the current 
author in \cite{RRV}.  We also discuss a new connection between an 
arithmetic condition of the additive submonoid $\langle 1/\lam_{1}, 
\ldots , 1/\lam_{n} \rangle$ called almost quasinormality and the 
regularity of the Rees algebra $R[It]$ in codimension one.

\section{Preliminaries}\label{sec:prelim}
\subsection{Integral Closure of Graded 
Rings}\label{subsec:intclosgrrings}

In this subsection we present some basic material about the integral 
closure of graded rings, which we will need when we talk about Rees 
algebras.  The standard technique for deducing results
about an ideal $I$  in a  ring $R$ 
is to pass to the Rees algebra $R[It]$ of $I$, which is the 
$R$-subalgebra of the polynomial algebra $R[t]$ generated by $\{ at 
\mid a \in I \}$.  We set the stage by proving a standard result about 
 the integral closure of a graded ring.  
 
 \smallskip
\textbf{Conventions.} All rings are assumed to be commutative with
identity and all ring homomorphisms preserve the identity.  
We let $\Z_+$ denote the set of positive integers, 
$\N$ the set of nonnegative integers,
$\Q_\ge$ the set of nonnegative rational numbers, $\Q_+$ the set of 
positive rational numbers, 
$\R_{\ge}$ the set of nonnegative real numbers, and $\mathbf{e}_1, 
\dots,
\mathbf{e}_n$ the standard basis vectors in $\R^n$.  We write 
$\balpha 
\le_{pr} \bbeta$ for vectors $\balpha = (a_{1},\ldots,a_{n}), \bbeta 
= 
(b_{1},\ldots,b_{n}) \in \R^{n}$ provided that $a_{i} \le b_{i}$ for 
 $1 \le i \le n$.  Thus $\balpha <_{pr} \bbeta$ means that $a_i \le 
b_i$ 
for all $1
\le i \le n$ and $a_j < b_j$ for some $1 \le j \le n$.  
By a semigroup we mean a set with an associative binary
operation.  By a monoid we mean a semigroup with identity.  We say an
abelian monoid $G$ is torsion-free provided that whenever
$f \ne g$ are elements of $G$ then $nf
\ne ng$ for all $n \in \Z_+$.  By a totally ordered abelian 
 monoid we mean an abelian monoid that is totally ordered and has the 
 property that whenever $f < g $ and $h$ are elements of $G$ then 
 $f+h < g + h$. For elements $g_1	, \ldots, g_r$  of an abelian monoid $G$
we let $\langle g_1, \dots, g_r \rangle$ denote
 the submonoid of $G$ generated by $g_1	, \ldots, g_r$.  If $G$ is an
additive abelian monoid and $K$ is any ring, we think of the monoid
ring
$K[G] = \oplus_{g \in G} Kx^g$ as a $G$-graded ring and say that the
homogeneous element $x^g$ of $K[G]$ has exponent $g$. If $G$ is a totally ordered
abelian monoid then $H=\grp(G)$ is naturally a totally ordered abelian
group.  If
$R$ is a $G$-graded ring we may regard $R$ as an $H$-graded ring by setting $R_g =
\{0\}$ for $g \in H \setminus G$.  By the integral closure, or
normalization,
$\overline{R}$ of a reduced ring
$R$ we mean the integral closure, or normalization, of $R$ in its total 
quotient ring.

Our first result is about the integral closure of a reduced graded ring
in a reduced graded extension ring.
 
\begin{theorem}\label{thm:intclosure}
 Let $G$ be a totally ordered abelian monoid and $R 
 \subseteq S$ be an extension of reduced $G$-graded rings.  Let $U$
denote 
 the set of homogeneous $S$-regular elements of $R$ and assume that 
 $U^{-1}R$ is integrally closed in $U^{-1}S$.  Then, the 
 integral closure of $R$ in $S$ is again $G$-graded.

\end{theorem}

\begin{proof}
 First let us assume that $R$ is Noetherian.  Suppose that 
 $s \in S$ is integral over $R$. For any element $s \in S$ let 
 in$(s)$ denote the initial component of $s$, that is the (nonzero) 
 homogeneous component of least degree.  Then $R[s]$ is a finite 
 $R$-submodule of  $U^{-1}R$.  Hence there exists an element $u \in U$ 
 such that $us^{i} \in R$ for all $i \ge 0$.  Since $S$ is reduced and
$u$ is not a zero divisor on $S$, the initial 
 component of  $us^{i}$ is $u$in$(s)^{i}$ and we must have $u$in$(s)^{i} 
 \in R$ for all $i \ge 0$.  Thus $R[\init(s)] \subseteq 
 Ru^{-1}$.  Since $R$ is Noetherian we may deduce that 
 $R[\init(s)]$ is a finitely-generated $R$-module and hence $\init(s)$ 
 is integral over $R$.  Hence $s - \init(s)$ is integral over $R$ and
by induction on the number of (nonzero) homogeneous components we may 
assume that each homogeneous component of $s - \init(s)$ is integral 
over $R$.

 Now consider the general case.  Suppose $s \in S$ satisfies the 
 equation
 $$ s^{n}+r_{1}s^{n-1}+ \cdots + r_{n} = 0 .$$
 As in the Noetherian case, there exists an element $u \in U$ 
 such that $us^{i} \in R$ for all $i \ge 0$ and hence $u\init(s)^{i} 
 \in R$ for all $i \ge 0$.  Let $R'$ be the $\Z$-subalgebra of $R$ 
 generated by $u$, the elements $us_{g}$, where $s_{g}$ is a 
 homogeneous component of $s$, and the homogeneous components of the 
 coefficents $r_{i}$ of the above equation.  By the proof of the 
 Noetherian case, $\init(s)$ is integral over $R'$ and by induction 
 on the number of homogeneous components every homogeneous component of 
 $s$ is integral over $R'$, hence is integral over $R$.
\end{proof}

\begin{rmk}\label{rem:intcl1}
    We point out that the hypotheses of the above theorem are 
    satisfied if $G$ is a finitely-generated totally ordered abelian
    monoid and
    $R \subseteq S$ is an extension of 
    $G$-graded integral domains with $S$ contained in the quotient field of
				$R$.  
    In this case we have $H := \grp(\balpha \in G \mid R_{\balpha} \ne 0
)$ is
				a finitely generated, torsion-free abelian group and hence is
				isomorphic to $\Z^d$ for some positive integer $d$.  The
    multiplicative subset $U$ above consists
				of the nonzero 
    homogeneous elements of $R$.   We have $U^{-1}R \cong K[H] \cong K[x_{1}, 
    \ldots, x_{d}, x_{1}^{-1}, \ldots, x_{d}^{-1}]$ is a ring of 
    Laurent polynomials over a field by \cite[Theorem 1.1.4]{GoWat}. With 
    this we can prove an absolute result about the integral closure 
    of a graded integral domain.  The corresponding assertion for 
    $\Z$-graded rings appeared in \cite[Lemma 2.1]{LV1}.
\end{rmk}

 \begin{cor}\label{cor:intcl}
     Let $G$ be a totally ordered abelian monoid and $R$ a $G$-graded 
     reduced ring with finitely many minimal primes.  Then, the integral 
     closure of $R$ is again $G$-graded.
 \end{cor}
 
 \begin{proof}
     First assume that $R$ is an integral domain. Let $U$ denote the set 
     of nonzero homogeneous elements of $R$ and let $H = \grp(G)$.  
     We will regard $R \subseteq U^{-1}R$ as an extension of $H$-graded
     rings.  
     If we can show that $U^{-1}R$ is integrally closed, then it will 
     follow that the integral closure of $R$ is the integral closure 
     of $R$ in $U^{-1}R$.  Then we can apply Theorem \ref{thm:intclosure} to 
     the extension $R \subseteq U^{-1}R$ to deduce the conclusion of 
     this corollary.

     Let $K$ denote the quotient field of $R$ and suppose that 
     $r/s \in K$ is integral over $U^{-1}R$.  Let $H'$ be the 
     subgroup of $H$ generated by exponents of the homogeneous
     components of $r, s$, 
     and the numerators and denominators of an equation of integral 
     dependence.  Let
     $A = \oplus_{g\in H'}
					R_{g}$ and $V$ be the set of 
     nonzero homogeneous elements of $A$. View $A \subseteq
					V^{-1}A$ as an extension of $H'$-graded domains. Then
					$r/s$ is in 
     the integral closure of $A$ and hence is in $V^{-1}A$ by 
     \cite[Theorem 1.1.4]{GoWat}.  Hence $r/s \in U^{-1}R$ as 
     asserted.

     In the general case, if $P_{1}, \ldots P_{\ell}$ are the minimal 
     primes of $R$, then each $P_{i}$ is homogeneous and 
     $\overline{R} = \overline{R/P_{1}} \times \cdots \times  
     \overline{R/P_{\ell}}$.  Since each factor $\overline{R/P_{i}}$ 
     is again $G$-graded so is $\overline{R}$.     
  \end{proof} 
 
  \medskip
    We will now state a version of Theorem \ref{thm:intclosure} that
				doesn't require that 
    the rings be integral domains or even reduced.

\begin{theorem}\label{thm:intclosureHS}
    Let $G$ be a cancellative, torsion-free abelian 
    monoid and $R \subset S$  an extension of $G$-graded rings.
    Then, the integral closure of $R$ in $S$ is again $G$-graded.
\end{theorem}

\begin{proof}
    This is a generalization due to Y. Yao \cite [Theorem 
    1.7.5]{HunSw} of the corresponding assertion when $G$ 
    is $\Z^{n}$ \cite [Theorem 1.7.3]{HunSw}.  Essentially the result 
    for $\Z^{n}$-graded rings is proved by a clever use of Vandermonde 
    matrices and then proved in the more general setting by reducing 
    to the $\Z^{n}$-case.   
\end{proof}               

\medskip

\subsection{Integral Closure of Homogeneous 
Ideals}\label{subsec:intclgradedids}

We now recall some fundamental results on the integral closure of an
ideal $I$ in a ring $R$ and the integral closure of a homogeneous ideal
in a graded ring.

\begin{defn}  Let $I$ be an ideal in a  ring $R$ and $x \in R$.  We 
say that $x$ is integral over $I$ if there is a positive integer $n$ 
and elements $a_{j} \in I^{j} (j=1, \ldots, n)$ such that
\[ x^{n}+a_{1}x^{n-1} + \cdots + a_{n} = 0. \]  The integral 
closure  $\overline{I}$ of $I$ is defined to be the set of all 
elements of $R$ that are integral over $I$.  The ideal $I$ is said to 
be integrally closed if $I = \overline{I}$.  The ideal $I$ is said to 
be normal if every positive power of $I$ is integrally closed.
\end{defn}

For an ideal $I$ in a  ring $R$ and an element $x \in R$, $x$ is integral 
over $I$ if and only if $xt \in R[t]$ is integral over 
the Rees algebra $R[It]$, where  $R[It]$ is the 
$R$-subalgebra of the polynomial algebra $R[t]$ generated by $\{ at 
\mid  a \in I \}$.     We can regard $R[t]$ as a $\N$-graded ring, where 
$R[t]_{n}=Rt^{n}$ for each $n \in \N$.  By 
 Theorem \ref{thm:intclosure} the integral closure of $R[It]$ in 
 $R[t]$ is again $\N$-graded.  This point of view yields the following.
 
 \begin{proposition}\label{prop:intclosureid}
     Let $I$ be an ideal of a ring $R$.  Then, the integral closure 
     $\overline{I}$ of $I$ is again an ideal of $R$.  Furthermore, 
     $\overline{I}$ is integrally closed.
 \end{proposition}
 
 \medskip
 The above conclusion is also a consequence of the following 
 well-known theorem (e.g., see \cite{Rib} or \cite[Proposition 5.1.7
]{HunSw}).
 
 \begin{theorem}\label{thm:intclRees}
     Let $R$ be a ring and $t$ an indeterminate.  For an ideal $I$ 
     of $R$, the integral closure of $R[It]$ in $R[t]$ is the graded 
     ring
    $\bigoplus_{n \ge 0} \overline{I^{n}}t^{n}.$
 \end{theorem}
 
 If we are working in a $G$-graded ring $R$ and have a homogeneous ideal 
 $I$ we can prove that $\overline{I}$ is again homogeneous by 
 considering $R[t]$ to be a $G \times \N$-graded ring, where 
 $R[t]_{(g,n)}=R_{g}t^{n}$.  We now state the exact result.
 
 \begin{theorem}\label{thm:gridealintcl}
     Let $G$ be a cancellative, torsion-free abelian monoid and let 
     $I$ be a homogeneous ideal in the $G$-graded ring $R$.  Then, 
     $\overline{I}$ is again homogeneous.  
 \end{theorem}

 \begin{proof}
     Note that $G \times \N$ is again a torsion-free abelian monoid 
     so by Theorem \ref{thm:intclosureHS} the integral closure of $R[It]$ in 
     $R[t]$ is again $G \times \N$-graded. Since an element $x \in R$ is
     in $\overline{I}$ if and only if $xt \in  
     \sum_{g \in G} \overline{R[It]}_{(g,1)}$, we may conclude that
     $\overline{I}$ is $G$-homogeneous.    
 \end{proof}

 \section{Monomial Ideals}
 
 \subsection{Basics on Monomial Ideals}\label{subsec:basicsmonid}
 
 In this section we recall some preliminary results about  monomial 
ideals in a polynomial ring over a field of arbitrary characteristic.
 We give proofs of some well-known results for 
 expository purposes and because we do not know of a 
 complete and correct reference for all of these assertions.

\begin{notn}  Throughout this section
$R$  will denote the polynomial ring $K[x_1,\dots, x_n]$
 over a field $K$ and 
$\mathfrak{m}=(\xvec{n})$ will denote  the maximal homogeneous ideal of
$R$.  In this context, for a vector 
$\balpha = (a_1, \ldots, a_n) \in \N^n$ we let $x^{\balpha}$ denote the monomial
$x_1^{a_1} \cdots x_n^{a_n}$.  
\end{notn}

After some preliminaries on the 
integral closure of a monomial ideal we turn our attention to the 
question of normality of a monomial ideal.
We will consider the polynomial ring as an 
$\N^{n}$-graded ring where for a vector $\balpha \in \N^{n}$ we let 
$R_{\balpha} = Kx^{\balpha}$.  With this grading 
homogeneous ideals and monomial ideals coincide.

\begin{corollary}\label{cor:monideal}
     If $I 
    \subseteq K[\xvec{n}]$ is a monomial ideal in a polynomial ring 
    in $n$ indeterminates over a field $k$, then $\overline{I}$ 
     is again a monomial ideal.  Indeed, 
$$\overline{I}=( x^{\balpha} 
     \mid  x^{m\balpha} \in I^{m} \; \exists m  \in \Z_+ ).$$
 \end{corollary}
 \begin{proof}
     Let $R=K[\xvec{n}]$ and regard $R$ as an $\N^{n}$-graded ring 
     as above.     Then 
     $\overline{I}$ is again a monomial ideal by 
     Theorem \ref{thm:gridealintcl}.  Suppose a monomial $x^{\balpha}$ 
     is integral over $I$.  Consider an equation of integral dependence
   $$x^{n\balpha} + a_{1}x^{(n-1)\balpha}+ \cdots + a_{n}=0,$$
 where $a_{j} \in I^{j}$.  Replacing each $a_{j}$ by its homogeneous 
 component of degree $j\balpha$ we may assume each $a_{j}
=c_{j}x^{j\balpha} \in 
 I^{j}$ where $c_{j} \in K$.  We must have $c_{j} \ne 0$ for some $j$ 
 between 1 and $n$
 and hence $x^{j\balpha} \in I^{j}$ for some positive integer $j$. 
 \end{proof}
 
 There is a convex-geometric description of the set of exponents of 
 monomials appearing in the integral closure of a monomial ideal that 
 is well-known but whose proof is difficult to find.  To our knowledge
the result first appeared for convergent power series rings in
\cite{LT}.  We present the proof for polynomial rings below.

For the reader's convenience we recall some definitions necessary to our
exposition. For the fundamentals on the convex geometry we are using the
reader can consult \cite{BH}, \cite{St} or \cite{Z}.  Although the
concepts of affine combination and affine dependence can be defined for
an arbitrary vector space, we restrict our attention to real vector
spaces for our convex-geometric definitions.

\begin{defn} Let $V$ be a real vector space. We say a vector $\balpha$ in 
$V$ is an \emph{affine combination} of the vectors $\balpha_1, \ldots,
\balpha_p$ in $V$ if $\balpha$ is a linear combination $\balpha =
c_1\balpha_1 + \cdots + c_p\balpha_p$ and $\sum_{i=1}^p c_i = 1$.   
If, in addition, each $c_i \ge 0$ we say $\balpha$ is a \emph{convex
combination} of the vectors $\balpha_1, \ldots, \balpha_p$.  A nonempty
subset of $\R^n$ that is closed under affine combinations is called an
\emph{affine subspace}.  Analogously, a nonempty subset of $\R^n$ is said to
be \emph{convex} if it is closed under convex combinations.
\end{defn}

Notice that if $W$ is an affine subspace of the real vector space $V$ and
$\balpha \in W$ then the set $U := W - \balpha = \{ \bbeta - \balpha \mid 
\bbeta \in W \}$ is a linear subspace of $V$  and $W = \balpha + U$.  The
subspace $U$ is uniquely determined by $W$ and by the dimension
$\dim(W)$ of $W$ we mean the dimension of $U$.

\begin{defn}  Let $S$ be a subset of a real vector space $V$. The set of all
affine combinations of vectors in $S$
is called the \emph{affine hull} $\aff(S)$ of $S$.  The set of all
convex combinations of vectors in $S$ is
called the \emph{convex hull} $\conv(S)$ of $S$. The set of all
nonnegative linear combinations (i.e., the coefficients are
nonnegative) of vectors in $S$ is called the \emph{positive cone}
$\pos(S)$ of
$S$.
\end{defn}

Observe that the affine hull of $S$ is the smallest affine
subspace containing $S$.  Similarly, the convex hull of $S$ is the
smallest convex subset containing $S$.

\begin{defn} We say a subset $\{ \balpha_1, \ldots, \balpha_p \}$ of a
 real vector space $V$ is \emph{affinely dependent} if there exist scalars
$c_1,
\ldots c_p$, not all zero, such that $\sum_{i=1}^p c_i\balpha_i =
\mathbf{0}$ and $\sum_{i=1}^p c_i = 0$.  If $\{ \balpha_1, \ldots,
\balpha_p \}$ is not affinely dependent we say the set is \emph{affinely
independent}.
\end{defn}

\begin{rmk} Suppose that $\sum_{i=1}^p c_i\balpha_i =
\mathbf{0}$  is an equation of affine dependence and that $c_1 \ne 0$. 
After scaling by a factor of $1/c_1$ we may may and shall assume that
$c_1=1$.  Hence $\balpha_1 = -c_2\balpha_2 - \cdots - c_p\balpha_p$
expresses $\balpha_1$ as an affine combination of $\balpha_2,
\ldots, \balpha_p$.  Thus the set $\{\balpha_1, \ldots, \balpha_p \}$ is
affinely dependent if and only if some vector in the set is an affine
combination of the remaining vectors.  Also notice that
$\{\balpha_1, \ldots, \balpha_p \}$ is affinely independent 
if and only if   $\{ \balpha_2 - \balpha_1, \ldots,
\balpha_p - \balpha_1 \}$ is linearly independent.
\end{rmk}

 \begin{defn}\label{gamma} Let $X$ be any subset of 
$R=K[x_1,\ldots,x_n]$.
Then set
$$\Gamma(X) \; = \; \{\balpha \in \N{}^n \hspace{.1cm} \mid  
x^{\balpha} \in X\}.$$
\end{defn}

We refer to $\Gamma(X)$ as the {\em exponent set of $X$}. 
If $I$ is a monomial ideal then $\Gamma(I)$ is an ideal of the 
monoid $\mathbb{N}^n$ \cite[page 3]{Gil}. If $A$ is a subalgebra of
$R$ generated by monomials then $\Gamma(A)$ is a submonoid of
$\mathbb{N}^n$, and $A$ is isomorphic to the monoid ring 
$K[\Gamma(A)]$.

Recall that if $I$ is a monomial ideal, then its exponent set $\Gamma(I)$
is a subset of $\N^n$.  We also regard it as a subset of $\R^n$, when
this is convenient.

\begin{defn} 
 For an arbitrary subset $\Lambda$ of  $\R{}^n$ and a positive 
 integer $m$ we let  
\begin{eqnarray*}
m\Lambda &=&  \{\lambda_1 + \cdots + \lambda_m \hspace{.1cm} 
\mid  
\lambda_i \in
\Lambda \hspace{.2cm} (i=1, \dots , m) \}.
\end{eqnarray*}

If $\Lambda = \Gamma(I)$ (respectively, $\Gamma(A)$) then 
$\conv(\Lambda)$ will be denoted NP$(I)$ (respectively, NP$(A)$),
and will be referred to as the {\em Newton polyhedron of $I$} 
(respectively, of $A$). 
\end{defn}

A polyhedron may be defined as the intersection
of finitely many half-spaces in $\R^n$ for some positive integer $n$.  A
polyhedron can also be thought of as the sum of the convex hull of a finite
point set (a polytope) and the positive cone generated by a finite set of
vectors.  Indeed these are equivalent notions (e.g., see \cite[Theorem
1.2]{Z}).  An affine hyperplane in $\R^n$ divides $\R^n$ into two
closed half-spaces.  If
$\Sigma \subseteq \R^n$ is a polyhedron and
$H$ is an affine hyperplane in $\R^n$ such that
$\Sigma \cap H \ne \emptyset$ and $\Sigma$ is contained in one of these
half-spaces, then
$H$ is called a supporting hyperplane of $\Sigma$.  A face of $\Sigma$
is the intersection of $\Sigma$ and a supporting
hyperplane.  If $F$ is a
face and $\dim(\aff(F)) =
\dim(\aff(\Sigma)) - 1$ then $F$ is called a
facet of $\Sigma$.       

The Newton polyhedron of a monomial ideal is an example of an
unbounded polyhedron.  To prepare for the geometric characterization of the
integral closure of a monomial ideal we now state and prove
Carath\'{e}odory's Theorem.

\begin{theorem}[Carath\'{e}odory's Theorem]\label{thm:cara}  Let $V$ be a
finite-dimensional real vector space and let $S \subseteq
V$ be a subset.  Then, a vector $\balpha \in V$ is in $\conv(S)$ if
and only if there exist affinely independent vectors
$\{ \balpha_1, \ldots, \balpha_p \}$ in $S$ with $ \balpha  \in
\conv(\balpha_1, \ldots, \balpha_p)$.
\end{theorem}

\begin{proof}Say $\balpha = \sum_{i=1}^p	b_i \balpha_i $ is a convex
linear combination.   Without loss of generality we may and shall
assume that $\balpha_i \ne \mathbf{0}$ and
$b_i > 0$ for all $i$. Suppose that $\{ \balpha_1, \ldots, \balpha_p 
\}$ is affinely dependent.   We proceed
by induction on $p$, the case $p=1$ being vacuously true.  Since 
$\{ \balpha_1, \ldots, \balpha_p 
\}$ is affinely dependent there exist scalars $c_1, \ldots, c_p$, not
all zero, such that $
\sum_{i=1}^p c_i\balpha_i = \mathbf{0}$ and $ \sum_{i=1}^p c_i = 0$.  By
multiplying by $-1$ if necessary, we may and shall assume that $c_i >0$ for
some index $i$.  For each index
$i$ such that $c_i >0$ consider $c_i/b_i$ and choose an index $i^*$ such
that 
$c_{i^*}/b_{i^*} \ge c_i/b_i$ for all indices $i$ such that $c_i > 0$.
Notice that $c_{i^*}/b_{i^*}$ = max$\{c_1/b_1, \ldots, c_p/b_p\}$.
Since the coefficient $c_{i^*} \ne 0$ we have  $\balpha_{i^*} =
-\sum_{i
\ne i^*} (c_i/c_{i^*})\balpha_i$, where $\sum_{i \ne i^*} (c_i/c_{i^*}) =
-1$.  Hence
$$\balpha = \sum_{i \ne i^*}b_i\balpha_i - b_{i^*}\sum_{i \ne
i^*}(c_i/c_{i^*})\balpha_i 
=\sum_{i \ne i^*} \left(b_i -
b_{i^*}(c_i/c_{i^*})\right)\balpha_i,$$ 
where 
$$\sum_{i \ne i^{*}} \left(b_i -
b_{i^*}(c_i/c_{i^*})\right) = \sum_{i \ne i^*} b_i - b_{i^*}\sum_{i
\ne i^*}(c_i/c_{i^*}) \\
= \left(\sum_{i \ne i^*} b_i\right) + b_{i^*} = 1,$$

\noindent and $b_i -
b_{i^*}(c_i/c_{i^*}) \ge 0$ for $i \ne i^*$.  By induction, $\sum_{i \ne i^*} 
(b_i - b_{i^*}(c_i/c_{i^*}))\balpha_i$ is a convex combination of
affinely independent vectors in $S$. 
\end{proof}

We can now state and prove the geometric characterization of the 
integral closure of a monomial ideal.

\begin{theorem}\label{normalization}
{\rm (a) } Let $I$ be a monomial ideal in $R=K[x_1,\ldots,x_n]$. Then 
the
integral closure $\overline{I}$ of $I$ in  $R$ is 
the monomial ideal defined by
$\Gamma(\overline{I})=\NP(I)\cap\mathbb{N}^n$ (so that 
$\NP(I)=\NP(\overline {I}))$.  Furthermore $$\Gamma(\overline{I})=
\{ \balpha \in \N{}^n |   m\balpha \in m \Gamma(I) \mbox{ for 
some }
 m \in \Z_+ \}.$$

{\rm (b)}  Let $A$ be a subalgebra of $R$ generated by
a finite number of monomials. Then the
integral closure $\overline{A}$ of $A$ in  $R$ is 
the semigroup ring defined by 
$\Gamma(\overline{A})=\NP(A)\cap\mathbb{N}^n$.  Furthermore 
$\NP(A)$ is the positive cone spanned
by $\Gamma(A)$ (or by the exponents of a (finite) set of algebra 
generators of $A$) and $$\Gamma(\overline{A})=
\{ \balpha \in \N{}^n |   m\balpha \in \Gamma(A) \mbox{ for some }
 m \in \Z_+\}.$$
\end{theorem} 
 
\begin{proof}
{\rm (a)} Suppose first that $x^{\balpha} \in \overline{I}$.  Thus 
$x^{m\balpha} \in I^{m}$ for some positive integer $m$ by Corollary 
\ref{cor:monideal}.  Thus there exist monomials $x^{\balpha_{1}}, 
\ldots, x^{\balpha_{m}} \in I$ such that  $m\balpha = \balpha_{1} + 
\cdots + \balpha_{m}$.  Hence 
$$\balpha = \frac{1}{m}\balpha_{1}+ \cdots + \frac{1}{m}\balpha_{m} \in 
\NP(I) \cap \N^{n}.$$
 
Now assume that $\balpha \in \NP(I) \cap \N^{n}$.  By Carath\'{e}odory's 
Theorem there exist affinely independent vectors $\balpha_{0}, \ldots , 
\balpha_{m} \in \Gamma(I)$ such that $\balpha \in \conv (\balpha_{0}, \ldots
, \balpha_{m})$.  We may assume $\balpha = c_{0}\balpha_{0} + \cdots 
+ c_{m}\balpha_{m}$, where the coefficients are positive real numbers 
whose sum is 1.  Then 
$$\balpha - \balpha_{0} = c_{1}(\balpha_{1}-\balpha_{0}) + \cdots +
c_{m}(\balpha_{m}-\balpha_{0}),$$
has a solution in $\R^{m}$ and hence must have a solution in 
$\Q^{m}$.  Since the vectors \newline $\balpha_{1}-\balpha_{0}, \ldots, 
\balpha_{m}-\balpha_{0}$ are linearly independent the solution is 
unique.  Hence the coefficients $c_{j}$ are positive rational numbers.  
Letting $d$ be a common denominator for the coefficients and writing 
each $c_{j} = a_{j}/d$ with $a_{j} \in \Z_+$ we have $d\balpha = 
a_{0}\balpha_{0} + \cdots + a_{m}\balpha_{m}$ and hence 
$x^{d\balpha} \in I^{d}$ since $\sum a_j = d$.  

{\rm (b)} This is \cite[Proposition 6.1.2]{BH}. See also
\cite[3.1]{RR} for a form closer to what we want here.   
\end{proof}

With this geometric characterization and the algebraic 
characterization of Corollary \ref{cor:monideal} in hand, 
we have a pretty good understanding of the integral closure of a monomial ideal. 
A vector $\balpha \in \N^{n}$ is in $\Gamma(\overline{I})$ if and only if 
$m\balpha \in m\Gamma(I)$ for some positive integer $m$. 
Thus the monomial ideal $I$ is integrally closed if and only if whenever 
$\balpha \in \N^{n}$ and $m\balpha \in m\Gamma(I)$ for some positive
integer $m$, the  vector $\balpha \in \Gamma(I)$.

\subsection{Normal $\mathfrak{m}$-primary monomial ideals}

We now turn our attention to the question of normality for monomial 
ideals.  Recall that an ideal $I$ is normal if all its positive powers 
are integrally closed.  Thus a monomial ideal $I$ in a polynomial  ring
$R =K[\xvec{n}]$ is normal if and only if the Rees algebra 
$R[It]$ is a normal domain by Theorem \ref{thm:intclRees}.  Although we
have a good characterization of an  integrally closed monomial ideal no
such characterization exists for a  normal monomial ideal, nor does it
appear that a general characterization can ever be given.   Thus we
will focus on a special class of monomial ideals, which we will now introduce. 
One can check
whether a particular monomial ideal is normal using \cite{normaliz}.

\begin{notn} Let $\blam = (\lam_1, \ldots, \lam_n)$ denote a
vector of positive integers.  
Set 
$J(\blam) =(x_{1}^{\lam_{1}},\ldots,x_{n}^{\lam_{n}})$ and 
$I(\blam)=\overline{J(\blam)}$.  
\end{notn}

After a general result on the 
normality of a monomial ideal we turn our attention to the 
question of normality of the $\mathfrak{m}$-primary ideal $I(\blam)$.

We begin by recalling some recent results of Reid, Roberts, and the 
current author in \cite{RRV}.

\begin{proposition} \cite[Proposition 3.1]{RRV}  \label{proples} Let $I
\subseteq R=K[\xvec{n}]$ be a  monomial 
ideal. If $I^{m}$ is integrally closed for $m = 1,\ldots,n-1$, then 
$I$ is 
normal.
\end{proposition}

Thus an integrally closed monomial ideal in $K[x_1,x_2]$ is automatically
normal (this is also a consequence of Zariski's work on complete ideals,
which is described in Appendix 5  of \cite{Zar}). Thus the first case of
interest is an integrally closed monomial ideal in
$K[x_1,x_2,x_3]$.  An integrally closed monomial ideal $I$ in
$K[x_1,x_2,x_3]$ is normal provided that $I^2$ is also integrally closed.  

\begin{notn} \label{notation4}
 Let $L= \lcm(\lam_{1},\ldots,\lam_{n})$, 
$\omega_{i}=L/\lam_{i}$, 
$1/\blambda=(1/\lambda_1,\ldots,1/\lambda_n)$ 
and $\bomega = 
(\omega_{1},\ldots,\omega_{n})$, so that $L/\blambda=\bomega$. 
We will denote $\Gamma(I(\blam))$ 
(Definition {\rm \ref{gamma})} simply by $\Gamma$.
\end{notn}

Observe that 
$\NP(I(\blam))=\NP(J(\blam))$ has one bounded facet 
with vertices $\lam_{1}\e_{1},\ldots,\lam_{n}\e_{n}$. 
For a general vector $\balpha=(a_1,\ldots,a_n)$ of $\R^n$ the hyperplane 
defined by the equation  
$ (1/\blambda) \cdot \balpha =1$
passes through these vertices, and upon multiplication by $L$, the
equation of  this hyperplane becomes $\bomega\cdot\balpha=L$.   
We point out that 
$\NP(I(\blam)^p)= \{ p\balpha \mid \balpha \in \NP(I(\blam)) \}$ also has
one bounded facet  with vertices $p\lam_{1}\e_{1},\ldots,p\lam_{n}\e_{n}$
and supporting hyperplane $\bomega\cdot\balpha=pL$ by \cite[Lemma 2.5]{RRV}.
For a precise  description of the faces of
the Newton polyhedron of a general monomial ideal $I$ the reader can consult
\cite{RV}. These remarks lead to the two lemmas that follow.

\begin{lemma} \cite[Lemma 4.1]{RRV} \label{notation4a} 
The exponent set $\Gamma$ of the monomial ideal $I(\blam)$ can be
described by $\Gamma =
\{ \balpha\in \N^n \mid  (1/\blambda) \cdot \balpha \ge 1\} =  
\{ \balpha \in \N^{n} \mid \bomega \cdot \balpha
\ge L\}$. 
\end{lemma}
The next lemma gives necessary and sufficient conditions for 
$I(\blam)$ to be normal in terms of the exponent set $\Gamma$. 

\begin{lemma} \cite[Lemma 4.3]{RRV} \label{lem reductions} For the ideal $I(\blam)$
defined  above the  following are equivalent.  
\begin{enumerate}
\item[{\rm (a)}] $I(\blam)$ is normal. 
\item[{\rm (b)}]   Whenever $\bomega \cdot \balpha \ge pL$ for 
$\balpha 
\in \N^{n}$ and $p \in \N$, there exist vectors \\ $\bbeta_{j} \in 
\Gamma \; 
(j=1,\ldots,p)$ 
 such that   $\balpha = \sum \bbeta_{j}.$
 \item[(c)]  For all integers $p$ with $1\leq p<n$ and vectors
$\balpha=(a_1,\ldots,a_n)\in\mathbb{N}^n$ with $\lambda_i>a_i  \; \; (i=1, \ldots,
n)$ satisfying $\bomega \cdot \balpha \ge pL$, there exist vectors $\bbeta_{j}
\in
\Gamma
\; (j=1,\ldots,p)$ 
 such that   $\balpha = \sum \bbeta_{j}.$
  
\end{enumerate}
\end{lemma}

Assigning deg($x_{i}) = \omega_{i} \; (i=1,\ldots, n)$ we now 
have $I(\blam) = R_{\ge \, L} := \bigoplus_{n \ge L} R_n$. Rings of this
type have recently been studied by K. Smith and her collaborators (e.g.,
see \cite{F}).  One fruitful consequence of this point of view is the
following result.

\begin{cor} \cite[Corollary 4.4]{RRV} Let $\blam=(\lam_{1},\ldots,\lam_{n}) \in
\Z_+^{n},\  n\geq 3,$ 
and suppose that $\gcd(\lam_{1},\ldots,\lam_{n}) >n-2$.  
Then the 
monomial ideal $I(\blam) \subseteq K[x_{1},\ldots,x_{n}]$ is normal. In
particular, if $n=3$ and the integers $\lam_{1},\lam_{2},\lam_{3}$ are not
relatively prime, then the ideal $I(\blam)$ is normal.
\end{cor}

In \cite{BG}  Bruns-Gubeladze define a  submonoid 
$S$ 
 of $\Q_{\ge}$ to be  \emph{1-normal} if  
whenever $x \in S$ and $x \leq p$ for some $p \in \N$, there exist 
rational numbers $y_{1},\ldots, y_{p}$ in $S$ with $y_i\leq 1$ for 
all $i$ 
such  that $x = y_{1} + \cdots + y_{p}$. Then they relate the 
normality of
$K[S(\blam)]$ to the 1-normality of the submonoid  
$\Lambda$  of $\mathbb{Q}_{\geq}$.  
In \cite{RRV} their program was modified as follows.

\begin{defn}
A submonoid  $S$  of $\Q_{\ge}$ is \emph{quasinormal} provided that 
whenever $x \in S$ and $x \ge p$ for some $p \in \N$, there exist 
rational numbers $y_{1},\ldots, y_{p}$ in $S$ with $y_i\geq 1$ for 
all $i$ 
such that $x = y_{1} + \cdots + y_{p}$.
\end{defn}
 
We have the following necessary condition for $I(\blam)$ to be normal.

\begin{lemma} \cite[Lemma 4.6, Proposition 4.7]{RRV}
\label{ifnormalthenquasinormal} Let 
$\Lambda = \langle 1/\lam_{1},\ldots,1/\lam_{n}\rangle$, the additive 
submonoid of $\Q_{\ge}$   generated by $1/\lam_{1},\ldots,1/\lam_{n}.$
If $I(\blambda)$ is normal then $\Lambda$ is quasinormal. Furthermore, if
integers $\lam_{1},\ldots,\lam_{n}$ are pairwise relatively prime the 
converse is true. 

\end{lemma}

Thus in the special case where the integers $\lam_{1},\ldots,\lam_{n}$ are
pairwise relatively prime, the normality condition  on the 
$n$-dimensional monoid $\Gamma$ is reduced to the quasinormality 
condition on the 1-dimensional monoid $\Lambda$.  Notice that in
$K[x_1,x_2,x_3]$  the only case where we do not have a good
answer to the question of when
$I(\blam)$ is normal is when $\gcd(\lam_1,\lam_2,\lam_3)=1$ and the
integers $\lam_1,\lam_2,\lam_3$ are not pairwise relatively prime.  This
remaining case is being investigated by the author and H. Coughlin.

The next result enables one to quickly conclude that certain monomial ideals
$I(\blam)$ are not normal and leads to another definition.

 \begin{proposition} \cite[Propositions 4.8 \& 4.9]{RRV}   Let $\blam \in
\Z_+^{n}$.
If $\Lambda = \langle 1/\lam_{1},\ldots,
 1/\lam_{n} \rangle$ is quasinormal, then $1 + 1/L \in \Lambda$ {\rm (}and hence, 
$L + 1 \in \langle 
 \omega_{1},\ldots,\omega_{n}\rangle$ {\rm )}.
 \end{proposition}

This condition leads us to our next definition.

\begin{defn}  We say that the submonoid  $\Lambda = 
\langle 1/\lam_{1},\ldots,
 1/\lam_{n} \rangle$  of $\Q_{\ge}$ is
\emph{almost quasinormal} provided that $1 + 1/L \in \Lambda$, equivalently, if 
$L + 1 \in \langle \omega_1, \ldots, \omega_n \rangle$.
\end{defn}

For the discussion below we need to recall yet another definition.

\begin{defn}  An \emph{affine semigroup} is a finitely generated monoid
that is isomorphic to a submonoid of $\Z^n$, for some positive integer
$n$. 
\end{defn}

At first almost quasinormality was viewed primarily as a way to produce
examples of integrally closed but not normal monomial ideals.  However, it
turns out that the condition is closely related to whether the Rees algebra of
the ideal satisfies condition ${\rm R}_1$ of Serre.  In order to discuss this
connection we observe that the Rees algebra $R[It]$ of a monomial ideal
$I$ of the polynomial ring $R = K[\xvec{n}]$  can always be
identified with an affine  semigroup ring over $K$.  Namely, if $I =
(x^{\bbeta_1}, \ldots, x^{\bbeta_r})$ and we let $S(I) = \langle (\e_1,0),
\ldots, (\e_n,0),(\bbeta_1,1), \ldots, (\bbeta_r,1) \rangle \subseteq
\N^{n+1}$, then $R[It] \cong K[S(I)]$.  If $I=I(\blam)$ then $S(I)$ may be
described as the submonoid of
$\N^{n+1}$ generated by
$$ \{ (a_1,\ldots,a_n,d) \in \mathbb{N}^{n+1} \mid
a_1/
\lambda_1 + \cdots + a_n/ \lambda_n \ge d \mbox{ for } d \le 1 \}.$$

Throughout the rest of this section $S$ will denote an affine semigroup,
				$K$ a field,
    $\mathcal{R}=K[S]$
    the affine semigroup ring over $K$, and $C=\pos(S)$ the positive cone 
    generated by $S$.  Each supporting hyperplane $H$ of $C$ is a
    linear subspance and is the zero set of some integral linear form
    $\sigma$. If the coefficients of $\sigma$ are relatively prime
    integers we refer to $\sigma$ as the primitive linear form defining
    $H$; this form is unique up to a factor of $-1$. Let
    $F$ be a facet of $C$ and $\sigma_{F}$ the 
    corresponding primitive linear form.  Let $H_{F}= \{ \balpha \in 
    \R S \mid \sigma_{F}(\balpha) = 0 \}$.  For a facet $F$ of $C$ let $P_F$
    denote the  height one  monomial prime generated by 
    the set  of monomials $\{ x^{\bbeta} \mid \bbeta \in S \setminus F  
    \} = \{ x^{\bbeta} \mid \bbeta \in S, \sigma_F(\bbeta) >0 
    \}$.  For definitions and details the reader
				should consult \cite[Chapter 6]{BH}.

The following two results were communicated to the author by W. Bruns.

\begin{lemma}\label{lem:bruns}
     Let $S$ be an affine semigroup and $F$ be a facet of the positive cone
				$C$ of $S$ such that $\mathcal{R}_P$ is regular, where $P = P_F$.
    Then we have:
\begin{itemize}
   \item[i.] there 
    exists an element  $\bbeta \in S \setminus F$ such that 
    $\sigma_{F}(\bbeta)=1$; and 

   \item[ii.] $\mathrm{grp}(S \cap F ) =\mathrm{grp}(S) 
\cap  H_{F} $.
\end{itemize}
\end{lemma}

\begin{proof}
     Let $P = P_F$.   By our assumption $\mathcal{R}_{P}$ is a discrete rank
				one valuation ring and  $P\mathcal{R}_{P}$ is generated by a single 
    monomial $x^{\bbeta}$.  Let $\mathrm{v}$ be the associated 
    valuation.   Then, $1= \mathrm{v}(x^{\bbeta})=\sigma_{F}(\bbeta)$.  
    Now let $\balpha, \, \balpha' \in S$ be such that $\balpha - \balpha'
    \in 
    F$. Then we must have
    $\balpha - \balpha' \in \mathrm{grp}(S) \cap H_{F}$.  Conversely 
    assume that $\balpha, \, \balpha' \in S$ and suppose that $\balpha - 
    \balpha' \in H_{F}$.  Then let $m = \sigma_{F}(\balpha) =  
    \sigma_{F}(\balpha')$.  If $m = 0$ then we have $\balpha, \, \balpha'
    \in 
    S \cap F$ and hence $\balpha - \balpha'  \in \mathrm{grp}(S \cap 
    F)$.  Suppose that $m > 0$.  Then $x^{\balpha}/x^{m\bbeta}$ 
    and  $x^{m\bbeta}/x^{\balpha'}$ are invertible monomials in 
    $\mathcal{R}_{P}$ since $P\mathcal{R}_{P}=x^{\bbeta}\mathcal{R}_P$. 
			 Hence $\balpha - \balpha' = (\balpha 
    - m\bbeta) + (m\bbeta - \balpha') \in \mathrm{grp}(S \cap F)$ as 
    asserted.
\end{proof} 
\begin{proposition} \label{prop:bruns}
    The monoid ring $K[S]$ of an affine semigroup $S$
     satisfies condition $\mathrm{R}_{1}$ of 
    Serre iff the following two conditions hold for every facet $F$ of the
    positive cone $C$ of $S$:
    \begin{itemize}
        \item [i.]  there exists an element $\bbeta \in S$ such that 
        $\sigma_{F}(\bbeta)=1$; and
    
        \item  [ii.]  $grp(S \cap F)  = \grp(S) \cap H_{F} $.
    \end{itemize}
\end{proposition}

\begin{proof}
    First assume that $\mathcal{R} = K[S]$ satisfies condition $\mathrm{R}_{1}$
    of Serre.  Let $P = P_F$ be a height one  prime corresponding to 
    the facet $F$ of $C$.  By Lemma \ref{lem:bruns} conditions (i)-(ii) are
satisfied.
    
    To prove the converse, assume that conditions (i) - (ii) hold.  Let 
    us first notice that $\mathcal{R}$ satisfies condition $\mathrm{R}_1$ if
    and only $\mathcal{R}_P$
    is regular for the the prime ideals $P=P_F$ corresponding to 
    facets $F$ of $C$. In fact, 
    if $Q$ is a height 1 prime ideal different from all the $P_F$, it cannot 
    contain a monomial (the minimal primes of a monomial ideal are generated by monomials 
    and therefore among the $P_F$). Suppose $S = \langle \bbeta_{1}, 
    \ldots , \bbeta_{r} \rangle$ and let $\bbeta = \bbeta_{1} + 
    \cdots + \bbeta_{r}$ and $f = x^{\bbeta}$.   Then 
    $\mathcal{R}[f^{-1}] = K[S,f^{-1}] =
    K[\mathrm{grp}(S)]$ as follows. If $\balpha = m_{1}\bbeta_{1}+ \cdots + 
    m_{r}\bbeta_{r} \in S$ then $m(\bbeta_{1} + \cdots + \bbeta_{r}) 
    - \balpha  \in S$, where $m = \max\{m_{i}\}$. Thus 
    $-\balpha \in \langle S, -\bbeta \rangle$.
    So $\mathcal{R}_Q$ is a localization of $K[\grp(S)]$ which is Laurent
    polynomial ring, and therefore regular.
    
    Now suppose that $P$ is a height one prime of $\mathcal{R}$ containing a 
    monomial $x^{\balpha}$.  Then $P$ is a minimal overprime of 
    $x^{\balpha}$ and hence is a height one monomial prime.  Thus 
    $P = P_{F} = \left( x^{\balpha} \mid \balpha \in S \setminus F 
    \right) $ for some facet $F$ of $C$.  We wish to see that
    $\mathcal{R}_{P}$ 
    is regular.  By condition (i) there is a vector $\bgamma \in S 
    \setminus F$ with $\sigma_{F}(\bgamma)=1$.  We claim that 
    $P\mathcal{R}_{P} = \left(x^{\bgamma} \right)$.  For suppose that 
    $x^{\balpha} \in P\mathcal{R}_{P}$.  Then $\balpha \in S \setminus F$ and 
    hence $\sigma_{F}(\balpha) = m >0$.  Thus $\balpha -m\bgamma \in 
    \grp(S) \cap H_{F} = \mathrm{grp}(S \cap F)$ by 
    condition (ii).  So we have $x^{\balpha} = 
    \left(x^{\balpha}/ x^{m\bgamma} \right)x^{m\bgamma} \in 
    x^{\bgamma}\mathcal{R}_{P}$ as desired.
    \end{proof}

We now wish to see what this means for the Rees algebra of an ideal 
$I = I(\blam)$.   As
we alluded to in the preceding paragraphs the counterpart of the Newton
polyhedron of a monomial ideal in the convex-geometric description of the affine
semigroup ring $K[S(I)]$ is the positive cone
$\pos(S)$ that is associated with the monoid 
$S=S(I)$. The facets of $\NP(I)$ are cut out by the supporting hyperplanes
$H_{\bomega,L}, H_1, \ldots, H_n$ where $H_{\bomega,L} = \{\balpha \in \R^n
\mid \bomega \cdot \balpha = L \}$ and $H_i$ is the coordinate hyperplane $H_i
=
\{\balpha \in \R^n \mid \e_i \cdot \balpha = 0 \}$.  The facets of
$C=\pos(S)$ are cut out by the supporting hyperplanes $H_{\sigma}, H_1,
\ldots, H_{n+1}$  where $\sigma(\balpha, a_{n+1})=\bomega \cdot \balpha
- La_{n+1}$ and $ H_1,
\ldots, H_{n+1}$ are the coordinate hyperplanes in $\R^{n+1}$.  Since the
height one monomial primes in $R[It]$ and $K[S(I)]$ can be identified we
have the following description of the height one monomial primes of
$R[It]$.

\begin{lemma}  For a monomial ideal $I=I(\blam)$ the height one
monomial primes of $R[It]$ are as follows:
\begin{eqnarray*}
P_i  &= & (x_i) + (x^{\bbeta_j}t \mid \e_i \le_{pr} \bbeta_j) \; for \; (i=1,
\ldots, n);\\
P_{n+1} &=& (x^{\bbeta_1}t, \ldots, x^{\bbeta_r}t); and  \\
 P_{\sigma} &= & (x_1, \ldots, x_n) + (x^{\bbeta_j}t \mid
\sigma(\bbeta_j,1) > 0).
\end{eqnarray*}
\end{lemma}

We can now characterize which vectors $\blam$ lead to Rees algebras that
satisfy condition $\mathrm{R}_1$ of Serre.  First we show that
condition (ii) of \ref{prop:bruns} is always satisfied by an affine
semigroup of the form $S(I(\blam))$.

\begin{lemma} \label{lem:groups} Let $\blam \in \Z_+^n$. $I=I(\blam)$, and let $S =
S(I) \subseteq \N^{n+1}$.  Then,
$$\grp(S \cap H_{\sigma}) = \grp(S) \cap H_{\sigma}. $$
\end{lemma}

\begin{proof}  Since $(\e_1,0), \ldots, (\e_n,0),
(\lam_1\e_1,1)$ are in $S$ the vectors
$(\e_1,0), \ldots, (\e_n,0), \newline(\mathbf{0},1)$ are in $\grp(S)$ and
hence
$\grp(S) = \Z^{n+1}$.  

It is clear that $\grp(S \cap H_{\sigma}) \subseteq \grp(S) \cap
H_{\sigma}$.  To prove the opposite containment suppose that
$(\balpha,a_{n+1}) \in \Z^{n+1}$ is on the hyperplane $H_{\sigma}$,
i.e., $ \sigma(\balpha,a_{n+1})=0$, where $\balpha = (a_1, \ldots, a_n)$.
Notice that $S \cap H_{\sigma}$ is generated by those generators
$(\bbeta_i,1)$ of $S$ for which
$\sigma(\bbeta_i,1)=0$.  Also observe that  $\{(\lam_1\e_1,1), \ldots,
(\lam_n\e_n,1) \} \subseteq \grp(S \cap H_{\sigma})$.  Write $a_i =
q_i\lam_i +r_i$ in $\Z$ with $0 \le r_i < \lam_i$ for $i=1, \ldots, n$.
Then $(\balpha,a_{n+1})-\sum_{i=1}^n
q_i(\lam_i\e_i,1)=(r_1,\ldots,r_n,d)$ where $d=a_{n+1}-\sum_{i=1}^n
q_i$.  Thus it suffices to show that
$(r_1,\ldots,r_n,d) \in \grp(S \cap H_{\sigma})$.  But
$(r_1,\ldots,r_n,d)$ is necessarily in $S$ and
$\sigma(r_1,\ldots,r_n,d) =0$ implies $(r_1,\ldots,r_n,d) =
(\bbeta_{i(1)},1) + \cdots + (\bbeta_{i(d)},1)$ where
$(\bbeta_{i(j)},1)
\in S \cap H_{\sigma} \; (j=1, \ldots , d)$; no generator of the form
$(\e_i,0)$ can occur since one such summand would force the entire sum
to have positive value under $\sigma$.  Thus $(r_1,\ldots,r_n,d) \in
\grp(S \cap H_{\sigma})$ as desired. 
\end{proof}

With this lemma in hand we can prove the connection between the 
conditions of quasinormality and regularity in codimension one.

\begin{proposition} For a monomial ideal $I=I(\blam)$  of $K[\xvec{n}]$ the
Rees algebra
$\mathcal{R} := R[It]$ satisfies condition $\mathrm{R}_1$ of Serre if and
only if the additive semigroup $\Lambda = \langle 1/\lam_1 , \ldots,
1/\lam_n \rangle$ is almost quasinormal. 
\end{proposition}

\begin{proof}
We first show that $\mathcal{R}_{P_i}$ is a 1-dimensional regular local ring
for $i=1, \ldots, n+1$ by showing that the unique maximal ideal is principal.

Let $S = S(\blam)$.  Suppose that $1 \le i \le n$.  We claim that $x_i$
generates the unique maximal ideal of $\mathcal{R}_{P_i}$.  Consider $
\bbeta_j$, where
$\e_i \le_{pr} \bbeta_j$. Write $\bbeta_j = \e_i + \balpha$, where $\balpha
\in \N^n$ and choose $m \in \{1, \ldots, \hat{i}, \ldots,  n \}$. Then
$$x^{\bbeta_j}t = (x^{\balpha}t)x_i =
\frac{x^{\balpha}}{x_m^{\lam_m}}(x_m^{\lam_m}t)x_i \in
x_i\mathcal{R}_{P_i}.$$   Hence $P_i\mathcal{R}_{P_i} =
x_i\mathcal{R}_{P_i}$.  Notice that for the linear form
$\sigma_i(a_1, \ldots, a_{n+1}) = a_i$ cutting out the facet $F_i =
S \cap H_i$ we have $\sigma_i(\e_i,0)=1$.

Now consider $P := P_{n+1} = (x^{\bbeta_1}t, \ldots, x^{\bbeta_r}t)$.  We
claim that $P\mathcal{R}_{P} = x^{\bbeta_1}t\mathcal{R}_{P}$. 
This follows since
$$x^{\bbeta_j}t=\frac{x^{\bbeta_j}}{x^{\bbeta_1}}x^{\bbeta_1}t \in
x^{\bbeta_1}t\mathcal{R}_{P}.$$
Notice that for the
linear form
$\sigma_{n+1}(a_1, \ldots, a_{n+1}) = a_{n+1}$ cutting out the facet $F_{n+1}
=\pos(S) \cap H_{n+1}$ we have $\sigma_{n+1}(\bbeta_1,1)=1$.  Notice that 
$P\mathcal{R}_{P} = x^{\bbeta_i}t\mathcal{R}_{P}$ for any $i$ between 1 and
$r$ by the same argument. 

Thus by Proposition \ref{prop:bruns} we know
that
$\mathcal{R}$ satisfies condition $\mathrm{R}_1$ of Serre if and only
$\mathcal{R}_{P_{\sigma}}$ is regular.  By Lemma \ref{lem:groups} this is the
case if and only if there exists a generator of 
$S$ with $\sigma$-value   equal to 1; the generator is either 
of the form $(\e_i,0)$ where  $\bomega \cdot \e_i =\omega_i = 1$ or
$(\bbeta_j,1)$ where  $\bomega \cdot \bbeta_j - L = 1$.  To finish the proof
observe that
 $\Lambda$ is almost quasinormal if and only if there exists $\balpha \in
\N^n$ with $\bomega \cdot \balpha = L + 1$. This happens if and only if for
some
$\bbeta_j$ and some $\bgamma \in \N^n$ we have $\sigma(\bbeta_j,1) +
\sigma(\bgamma,0) =1$; notice that $\bgamma$ is either $\mathbf{0}$ or
$\e_i$ for some $i$.  Thus $\Lambda$ is almost quasinormal if and only if 
$S$ has a generator with
$\sigma$-value equal to 1.
\end{proof}

In closing we mention one last remarkable result from \cite{RRV}.
Our notation continues as usual: $R=K[x_1,\ldots,x_n]$ for a field 
$K$, $\blambda=(\lambda_1,\ldots,\lambda_n)$ for arbitrary positive 
integers 
$\lambda_j$, and $L=\lcm(\lam_1, \ldots, \lam_n)$.  Define $\blam' = 
(\lam_1, \ldots, \lam_{i-1},
\lam_i + \ell, \lam_{i+1} \ldots,
\lam_n)$, where $ \ell = \lcm(\lam_1, \ldots, \widehat{\lam_i},
\ldots, \lam_n)$.  There is no loss of generality in taking
$i=n$, so that
$\ell=\lcm(\lambda_1,\ldots,\lambda_{n-1})$ and $\blambda^\prime= 
(\lambda_1,\ldots,\lambda_{n-1},\lambda_n+\ell)$.
We now state the result.
\smallskip

\begin{theorem} \cite[Theorem 5.1]{RRV} \label{congruence} If
$I(\blambda^{\prime})$ is  normal then 
$I(\blambda)$ is normal. If $\lambda_n \geq \ell$ and  
$I(\blambda)$ is normal so is  $I(\blambda^{\prime})$.
\end{theorem}

\end{document}